\documentclass[12pt]{amsart}

\usepackage{amsmath}
\usepackage{amsthm}
\usepackage{amssymb}
\usepackage{eucal}
\usepackage{times}
\usepackage{euler}
\usepackage{eufrak}

\usepackage[dvips]{graphicx}

\title[Gaussian fluctuations of representations of wreath products]%
{Gaussian fluctuations of representations\\ of wreath products}

\author{Piotr \'Sniady}

\address{Piotr \'Sniady \\ Institute of Mathematics \\
University of Wroclaw \\ pl.~Grunwaldzki~2/4 \\ 50-384 Wroclaw
\\ Poland}

\email{Piotr.Sniady@math.uni.wroc.pl}

\numberwithin{equation}{section} \numberwithin{figure}{section}

\theoremstyle{plain}
\newtheorem{lemma}{Lemma}[section]
\newtheorem{theorem}[lemma]{Theorem}
\newtheorem{theoremanddefinition}[lemma]{Theorem and Definition}
\newtheorem{proposition}[lemma]{Proposition}
\newtheorem{corollary}[lemma]{Corollary}

\theoremstyle{definition}

\theoremstyle{remark}

\newtheorem{example}[lemma]{Example}

\newcommand{\A}{{\mathfrak{A}}}

\newcommand{\E}{{\mathbb{E}}}
\newcommand{\C}{{\mathbb{C}}}
\newcommand{\R}{{\mathbb{R}}}
\newcommand{\Z}{{\mathbb{Z}}}
\newcommand{\M}{{\mathcal{M}}}
\newcommand{\N}{{\mathbb{N}}}

\newcommand{\iloscklas}{{|\hat{G}|}}

 \newcommand{\partialS}{\mathcal{P}S}

 \DeclareMathOperator{\Cov}{Cov}
 \DeclareMathOperator{\tr}{tr}
 \DeclareMathOperator{\Tr}{Tr}
 \DeclareMathOperator{\End}{End}

\begin{document}

\begin{abstract}
We study the asymptotics of the reducible representations of the
wreath products $G\wr S_q=G^q \rtimes S_q  $ for large $q$, where
$G$ is a fixed finite group and $S_q$ is the symmetric group in $q$
elements; in particular for $G=\Z/2\Z$ we recover the
hyperoctahedral groups. We decompose such a reducible representation
of $G\wr S_q$ as a sum of irreducible components (or, equivalently,
as a collection of tuples of Young diagrams) and we ask what is the
character of a randomly chosen component (or, what are the shapes of
Young diagrams in a randomly chosen tuple). Our main result is that
for a large class of representations the fluctuations of characters
(and fluctuations of the shape of the Young diagrams) are
asymptotically Gaussian. The considered class consists of the
representations for which the characters asymptotically almost
factorize and it includes, among others, the left regular
representation therefore we prove the analogue of Kerov's central
limit theorem for wreath products.
\end{abstract}

\maketitle

\section{Introduction}
\subsection{Formulation of the problem: asymptotics of representations}
One of the classical problems in the theory of asymptotic
combinatorics concerns the limit behavior of representations of
large groups. The above formulation of the problem is very general
and vague, so let us be more specific. An example of a question
which fits into the above category is the following one: sequences
$(G_q)$ and $(\rho_q)$ are given, where for each $q\in\N$ we know
that $G_q$ is a finite group and $\rho_q$ is a finite-dimensional
representation of $G_q$; we decompose $\rho_q$ as a sum of
irreducible components and we ask about asymptotic properties of a
random summand. Of course, if the irreducible representations of the
groups $G_q$ do not have a nice common structure then it is not
clear which \emph{asymptotic properties of a random summand} could
be studied and for this reason we should restrict our attention to
\emph{`nice'} groups $G_q$ for which the irreducible representations
can be nicely described in a uniform way.

In this article we are concerned with the case when $G_q=G \wr
S_q=G^q \rtimes S_q $ is the wreath product of a fixed finite group
$G$ by a symmetric group $S_q$. However, before we come to this
topic in Section \ref{sub:Hyperoc-groups} we will have a closer look
(in Section \ref{sub:Symmetr-groups} below) on a more developed case
of the symmetric groups.

\subsection{Asymptotics of representations of the symmetric groups}
\label{sub:Symmetr-groups} The simplest example of such
\emph{`nice'} groups is the sequence of the symmetric groups
$S_q$---in this case the irreducible representations of $S_q$ are
enumerated by the Young diagrams with $q$ boxes. Our problem
therefore asks about the statistical properties of some random Young
diagrams as their size tends to infinity. This problem was studied
in much detail: Logan and Shepp \cite{LoganShepp} and Vershik and
Kerov \cite{VershikKerov1977} studied the situation when $\rho_q$ is
the left-regular representation (which corresponds to, so called,
Plancherel measure on the irreducible representations of $S_q$) and
they proved that the corresponding Young diagrams concentrate (after
some simple geometric rescaling) around some limit shape. Biane
\cite{Biane1998,Biane2001approximate} generalized this result and he
proved the concentration of the shapes for a very large class of
representations of $S_q$; he also found a connection between the
asymptotic theory of the representations of the symmetric groups and
Voiculescu's free probability theory \cite{VoiculescuDykemaNica}.

The above results can be viewed as an analogue of the law of large
numbers, therefore it was very tempting to check if some kind of
central limit theorem could be true. The first result of this kind
was found by Kerov \cite{Kerov1993gaussian} (see also Ivanov and
Olshanski \cite{IvanovOlshanski2002}) who proved that for the
left-regular representations
the fluctuations of a random Young diagram around the limit shape
are asymptotically Gaussian. Hora
\cite{Hora2002noncommutativeaspect,Hora2003noncommutativeKerov}
found the non-commutative aspect of this result. In a previous paper
\cite{Sniady2005GaussuanFluctuationsofYoungdiagrams} we proved that
an analogue of the Kerov's central limit theorem holds true for a
large class of representations $\rho_q$ with, so called,
\emph{approximate factorization of characters}.


\subsection{The main result: asymptotics of representations of the wreath products}
\label{sub:Hyperoc-groups} It seems that our understanding of the
asymptotic properties of the symmetric groups is reaching a level of
maturity and for this reason the time has come to have a look on
some other classical series of finite groups and in this article we
will concentrate on wreath products $G \wr S_q$, where $G$ is a
fixed finite group.

Since $G \wr S_q= G^q \rtimes S_q$ is equal to a semidirect product
of a symmetric group $S_q$ therefore the structures of the groups $G
\wr S_q$ and $S_q$ and of their irreducible representations are very
closely related to each other; for this reason our strategy in this
article is to provide a setup in which the known results and methods
concerning $S_q$
\cite{Sniady2005GaussuanFluctuationsofYoungdiagrams} could be
directly applied to $G \wr S_q$. We will do it in Section
\ref{sec:Prelimi}.

Section \ref{sec:main-result:-represe-with-approxi-factori-charact}
contains the main result, namely that the methods which we used for
the study of the symmetric groups
\cite{Sniady2005GaussuanFluctuationsofYoungdiagrams} can be extended
to the wreath products. To be more explicit: we will introduce a
class of representations with \emph{approximate factorization of
characters}; this class will turn out to be big enough to contain a
lot of natural examples and will be closed under some natural
operations on representations such as induction, restriction, outer
and tensor product. We will prove that for the representations with
approximate factorization of characters an analogue of the Kerov's
theorem \cite{Kerov1993gaussian} holds true, in other words the
fluctuations of a randomly chosen irreducible component are
asymptotically Gaussian.




\section{Preliminaries}
\label{sec:Prelimi}

%
%

\subsection{Partial permutations and disjoint products}
\label{subsec:algebraofconjugacy}
\label{sub:Partial-permuta-disjoin-product} Ivanov and Kerov
\cite{IvanovKerov1999} defined a partial permutation as a pair
$\alpha=(\pi,A)$, where $A$ (called support of $\alpha$) is any
subset of $\{1,\dots,q\}$ and $\pi:\{1,\dots,q\}\rightarrow
\{1,\dots,q\}$ is a bijection which is equal to identity outside of
$A$. The natural product of partial permutations is given by
$$ (\pi_1,A_1)(\pi_2,A_2)=(\pi_1 \pi_2,A_1\cup A_2).$$
Partial permutations form a semigroup $\partialS_q$; in this article
we are interested also in the corresponding semigroup algebra
$\C(\partialS_q)$ which should be regarded as an analogue of the
permutation group algebra $\C(S_q)$ equipped with some additional
structure.


The algebra $\C(\partialS_q)$ can be also equipped with a different
product $\bullet$, called disjoint product, given on generators by
$$(\pi_1,A_1)\bullet(\pi_2,A_2)=\begin{cases}
(\pi_1 \pi_2,A_1\cup A_2) & \text{if } A_1\cap A_2=\emptyset, \\
0 & \text{otherwise}. \end{cases}$$

\subsection{Conjugacy classes for symmetric groups}

\label{subsec:definicjasigma}

Let integer numbers $k_1,\dots,k_m\geq 1$ be given. We define the
normalized conjugacy class indicator to be a central element in the
group algebra $\C(S_q)$ given by
\cite{KerovOlshanski1994,BianeCharacters,Sniady04AsymptoticsAndGenus,
Sniady2005GaussuanFluctuationsofYoungdiagrams}
\begin{equation}
\label{eq:definicjasigma} \Sigma_{k_1,\dots,k_m}=
\sum_{a} (a_{1,1},a_{1,2},\dots,a_{1,{k_1}}) \cdots
(a_{m,1},a_{m,2},\dots,a_{m,k_m}), \end{equation}
where the sum runs over all one--to--one functions
$$a:\big\{ \{r,s\}: 1\leq r\leq m, 1\leq s\leq {k_r}\big\}\rightarrow \{1,\dots,q\}$$ and
$(a_{1,1},a_{1,2},\dots,a_{1,{k_1}}) \cdots
(a_{m,1},a_{m,2},\dots,a_{m,k_m})$ denotes the product of disjoint
cycles. Of course, if $q<k_1+\cdots+k_m$ then the above sum runs
over the empty set and $\Sigma_{k_1,\dots,k_m}=0$.

In other words, if $k_1\geq \dots \geq k_m$ we consider a Young
diagram with the rows of the lengths $k_1,\dots,k_m$ and all ways of
filling it with the elements of the set $\{1,\dots,q\}$ in such a
way that no element appears more than once. Each such a filling can
be interpreted as a permutation when we treat the  rows of the Young
tableau as disjoint cycles. It follows that $\Sigma_{k_1,\dots,k_m}$
is a linear combination of permutations which in the cycle
decomposition have cycles of length $k_1,\dots,k_m$ (and,
additionally, $q-(k_1+\cdots+k_m)$ fix-points). Each such a
permutation appears in $\Sigma_{k_1,\dots,k_m}$ with some positive
integer multiplicity depending on the symmetry of the tuple
$k_1,\dots,k_m$.

It is natural to regard $\Sigma_{k_1,\dots,k_m}$ as an element of
$\C(\partialS_q)$, where we treat each summand on the
right-hand-side of \eqref{eq:definicjasigma} as a partial
permutation with the support equal to the set $\{a_{r,s}\}$.

\subsection{Algebra of conjugacy classes}
The sequence of the algebras $\C(\partialS_q)$ forms an inverse
system with the maps $\C(\partialS_{q+1})\rightarrow
\C(\partialS_q)$ equal to restrictions. We consider a vector space
spanned by the family of elements $\lim_{\leftarrow}
\Sigma_{k_1,\dots,k_m} \in \lim_{\leftarrow} \C(\partialS_q)$. It is
easy to check that both the natural product and the disjoint product
of two such elements is again in this form therefore the above
vector space is an algebra both when as the product we take the
natural product (we will denote this algebra by $\A$) and when as the
product we take the disjoint product (we will denote this algebra by
$\A^{\bullet}$).

%

\subsection{Wreath product and its representations}
For a finite group $G$ its wreath product with the symmetric group
$G\wr S_q$ is a semidirect product $G^q \rtimes S_q$, where $S_q$
acts on $G^q$ by permuting the factors.

We denote by $\hat{G}$ the set of the (equivalence classes of)
irreducible representations of the finite group $G$. For an
irreducible representation $\zeta\in \hat{G}$ we denote by
$p_{\zeta}\in\C(G)$ the corresponding minimal central projection.

Let a Young diagram $\lambda$ with $|\lambda|$ boxes be given and
let $\rho_{\lambda}:S_{|\lambda|}\rightarrow \End(W_{\lambda})$ be
the corresponding irreducible representation of the symmetric group.
For an irreducible representation $\zeta:G\rightarrow
\End(V_{\zeta})$ we consider the tensor power representation
$\zeta^{\otimes |\lambda|}$ of the group $G^{|\lambda|}$ acting on
the vector space
\begin{equation}
\label{eq:reprezentacja-produktowa}  (V_{\zeta})^{\otimes |\lambda|}
.
\end{equation}
Please notice that \eqref{eq:reprezentacja-produktowa} is also a
representation of $S_{|\lambda|}$, where the symmetric groups acts
by permuting the factors; in this way
\eqref{eq:reprezentacja-produktowa} is a representation of the
semidirect product $G^{|\lambda|}\rtimes S_{|\lambda|}=G\wr S_q$.

Also, $W_{\lambda}$ is a representation of $G^{|\lambda|}\rtimes
S_{|\lambda|}$, where $S_{|\lambda|}$ acts on this space by
$\rho_{\lambda}$ and $G^{|\lambda|}$ acts trivially as identity.

Thus, the tensor product
\begin{equation}
\label{eq:reprezentacja-produktowa2} W_{\lambda}\otimes
(V_{\zeta})^{\otimes |\lambda|} \end{equation}
is also a
representation of $G\wr S_q$. It will be the basic building block in
the construction of all irreducible representations of $G\wr S_q$.
The following well known result gives a complete classification of
the irreducible representations of the wreath products.
\begin{proposition}
\label{prop:reprezentacje}
Let $\Lambda:\hat{G}\rightarrow
\mathbb{Y}$ be a function on the set of irreducible representations
of $G$ valued in the set of the Young diagrams $\mathbb{Y}$ such
that $\sum_{\zeta\in\hat{G}} |\Lambda(\zeta)|=q$.

The representation
\begin{equation}
\label{eq:reprezentacje-produktow} \rho_{\Lambda} =\Big(
\bigotimes_{\zeta \in \hat{G} } W_{\Lambda(\zeta)} \otimes
V_{\zeta}^{\otimes |\Lambda(\zeta)|}\Big) \Big\uparrow^{G^q \rtimes
S_q}_{\prod_{\zeta\in \hat{G} } \left(G^{|\Lambda(\zeta) |}\rtimes
S_{|\Lambda (\zeta)|}\right)}
\end{equation} is irreducible; furthermore every irreducible representation of
$G^q\rtimes S_q$ is of this form.
\end{proposition}

\begin{example}
We consider the case when $G=\Z/2 \Z=\{0,1\}$. Then there are two
irreducible representations of $\Z/2 \Z$, namely
$\widehat{\Z/2 \Z}=\{\zeta_1,\zeta_2\}$, where $\zeta_1$ is the trivial
representation and $\zeta_2$ is the ``alternating'' representation
\begin{align*} \zeta_1(0)&=1, & \zeta_1(1)&=1, \\ \zeta_2(0)&=1, & \zeta_2(1)&=-1.
\end{align*}
The spaces $V_1$, $V_2$ on which these representations act are
one-dimensional.

In the case when $\zeta=\zeta_i$ the representation
\eqref{eq:reprezentacja-produktowa2} of $(\Z/2\Z)\wr S_{|\lambda|}$
takes a simpler form
$$ W_{\lambda}\otimes (V_i)^{\otimes |\lambda|}\cong W_{\lambda}; $$
on which the group $S_{|\lambda|}$ acts by $\rho_{\lambda}$ and the
elements $(z_1,\dots,z_{\lambda})\in (\Z/2\Z)^{|\lambda|}$ act
either trivially as identity (if $\zeta=\zeta_1$) or by multiplying
by $(-1)^{z_1+\cdots+z_{|\lambda|}}$ (if $\zeta=\zeta_2$).

The function $\Lambda:\widehat{\Z/2 \Z}\rightarrow \mathbb{Y}$ can be
identified with a pair of Young diagrams $\lambda_1,\lambda_2$,
where $\lambda_i=\Lambda(\zeta_i)$. The representation
\eqref{eq:reprezentacje-produktow} can be written more
explicitly as
$$\rho_{\lambda_1,\lambda_2} =\Big( W_{\lambda_1} \otimes W_{\lambda_2} \Big)
\Big\uparrow^{(\Z/2 \Z)^q \rtimes S_q}_{\left((\Z/2 \Z)^{|\lambda_1 |}\rtimes
S_{|\lambda_1|}\right)\times \left((\Z/2 \Z)^{|\lambda_2 |}\rtimes
S_{|\lambda_2|}\right)}.
$$
\end{example}

\subsection{Normalized trace}

For a matrix $x\in \M_n(\C)$ we define its normalized trace
$$ \tr x=\frac{1}{n} \Tr x,$$
where $\Tr x$ denotes the usual trace of $x$. In this way the
normalized trace of the identity matrix fulfills $\tr 1=1$.

\subsection{The main homomorphism}

The following Lemma is of great importance for this article: it
provides a homomorphism from the tensor product $\A^{\otimes
\iloscklas}$ to $\C(G\wr S_q)$; in this way every representation of
the wreath product $G\wr S_q$ defines a representation of the tensor
product $\A^{\otimes |\hat{G}|}$ and the questions concerning the
asymptotic behavior of the wreath products are reduced to the
corresponding questions concerning the symmetric groups.

For an irreducible representation $\zeta\in\hat{G}$, the
corresponding minimal central projection $p_{\zeta}\in\C(G)$ and a
partial permutation $(\pi,A)\in
\partialS_q$ we set
\begin{equation}
\label{eq:homomorfizm1} \phi_\zeta(\pi,A)=  (r_{1} \times \cdots
\times r_{q} ) \pi \in \C(G^q \rtimes S_q),
\end{equation} where
$$r_{m}=\begin{cases} p_{\zeta} & \text{if\/ } m\in A, \\ 1 &
\text{otherwise}.
\end{cases}
$$
\begin{lemma}
\label{lem:glowny-homomorfizm} \
\begin{enumerate} \item The map $\phi_{\zeta}$ extends to a homomorphism of algebras
$\phi_{\zeta}:\C(\partialS_q)\rightarrow \C(G^q\rtimes S_q)$, where we consider
$\partialS_q$ equipped with the natural product;
\item Let us fix for a moment some numbering of the irreducible representations
$\hat{G}=\{\zeta_1,\dots,\zeta_{\iloscklas}\}$.  The map
$\phi:\big(\C(\partialS_q)\big)^{\otimes |\hat{G}|}\rightarrow
\C(G^q\rtimes S_q)$ given on simple tensors by
$$ \phi\big( (\pi_1,A_1) \otimes \cdots \otimes
(\pi_\iloscklas,A_{\iloscklas}) \big)= \prod_{1\leq j\leq
\iloscklas} \phi_{\zeta_{j}} (\pi_j,A_j) $$ is a homomorphism of
algebras.
\item Let $\rho_{\Lambda}$ be an irreducible representation of
$G\wr S_q$, as presented in Proposition \ref{prop:reprezentacje}.
For any $a_1,\dots,a_\iloscklas\in\A$
\begin{equation} 
\label{eq:formula-na-charaktery}
\tr \rho_{\Lambda}\big( \phi(a_1\otimes \cdots \otimes
a_\iloscklas)\big)=
\prod_{1\leq j\leq \iloscklas} \tr \rho_{\Lambda(\zeta_j)}(a_j),
\end{equation}
where
$\rho_{\Lambda(\zeta_j)}$ is the irreducible representation of\/
$S_{|\Lambda(\zeta_j)|}$ corresponding to the Young diagram
$\Lambda(\zeta_j)$.
\end{enumerate}
\end{lemma}
\begin{proof}
If a permutation $\pi$ leaves a set $A$ invariant then the elements
$r_1 \times \cdots \times r_q$ and $\pi   \in \C(G^q \rtimes S_q)$
contributing to \eqref{eq:homomorfizm1} commute. Let
$(\pi^{(1)},A^{(1)}),(\pi^{(2)},A^{(2)})\in\partialS_q$ and
$(\pi^{(3)},A^{(3)}):=(\pi^{(1)},A^{(1)}) (\pi^{(2)},A^{(2)})$; it
follows that \begin{multline*} \phi_\zeta(\pi^{(1)},A^{(1)})
\phi_\zeta(\pi^{(2)},A^{(2)})=  \\ \pi^{(1)} (r^{(1)}_1 \times
\cdots \times r^{(1)}_q) (r^{(2)}_1 \times \cdots \times r^{(2)}_q)
\pi^{(2)}= \\
\pi^{(1)} (r^{(3)}_1 \times \cdots \times r^{(3)}_q)  \pi^{(2)}=
\pi^{(1)}  \pi^{(2)} (r^{(3)}_1 \times \cdots \times r^{(3)}_q) =\\
\phi_\zeta(\pi^{(3)},A^{(3)})
 \end{multline*}
which finishes the proof of point (1).

In order to prove point (2) it is enough to prove that if
$\zeta^{(1)}\neq \zeta^{(2)}$ then the elements
$\phi_{\zeta^{(1)}}(\pi^{(1)},A^{(1)})$ and
$\phi_{\zeta^{(2)}}(\pi^{(2)},A^{(2)})$ commute. This holds true
because when the sets $A^{(1)},A^{(2)}$ are disjoint then
permutations $\pi^{(1)},\pi^{(2)}$ commute; when the sets
$A^{(1)},A^{(2)}$ are not disjoint then
$$\phi_{\zeta^{(1)}}(\pi^{(1)},A^{(1)})\
\phi_{\zeta^{(2)}}(\pi^{(2)},A^{(2)})=0
=\phi_{\zeta^{(2)}}(\pi^{(2)},A^{(2)})\
\phi_{\zeta^{(1)}}(\pi^{(1)},A^{(1)}).$$

Since $a_1,\dots,a_\iloscklas\in\A$ therefore $\phi(a_1\otimes \cdots \otimes
a_\iloscklas)\in \C(G^q\rtimes S_q)$ is central and Frobenius reciprocity can
be applied to calculate the characters of the induced representation
\eqref{eq:reprezentacje-produktow}. This shows that the left-hand side of
\eqref{eq:formula-na-charaktery} is equal to the character of the representation
$$\bigotimes_{1\leq j\leq \iloscklas} W_{\Lambda(\zeta_j)} \otimes
V_{\zeta_j}^{\otimes |\Lambda(\zeta_j)|}
$$
of the group
$$G_{\Lambda}:=\prod_{1\leq j\leq\iloscklas } \left(G^{|\Lambda(\zeta_j)
|}\rtimes
S_{|\Lambda (\zeta_j)|} \right)$$
evaluated on the element 
$ \phi(a_1\otimes \cdots \otimes
a_\iloscklas)\big\downarrow_{G_{\Lambda}}$.
Notice that the element of the form
$$ \phi\big( (\pi_1,A_1)\otimes \cdots \otimes
(\pi_{\iloscklas},A_{\iloscklas}) \big)\big\downarrow_{G_{\Lambda}}
$$ represents on this space as
$$ \begin{cases} 
\bigotimes_{1\leq j\leq \iloscklas}
\left[ \big(\Lambda(\zeta_j)\big) (\pi_j) \otimes 1^{\otimes
|\Lambda(\zeta_j)|}\right] 
& \text{if $A_j$ is a
subset of the set of the} \\ 
& \text{elements permuted by subgroup} \\
& \text{$S_{|\Lambda (\zeta_j)|}$ of $G_{\Lambda}$ for each $j$}; \\
0 & \text{otherwise.}
   \end{cases}
$$
It follows that the element $ \phi(a_1\otimes \cdots \otimes
a_\iloscklas)\big\downarrow_{G_{\Lambda}}$ represents as 
$$ \bigotimes_{1\leq j\leq \iloscklas}
\left[ \big(\Lambda(\zeta_j)\big) (a_j) \otimes 1^{\otimes
|\Lambda(\zeta_j)|}\right]  $$
which finishes the proof of point (3).
\end{proof}

We define the natural product on $\A^{\otimes \iloscklas}$ by
\begin{equation} \label{eq:iloczyn-naturalny-na-tensorach}
(a_1\otimes \cdots \otimes a_{\iloscklas}) (b_1\otimes \cdots
\otimes b_{\iloscklas})= a_1 b_1\otimes \cdots \otimes
a_{\iloscklas} b_{\iloscklas};\end{equation}  and the disjoint
product $\bullet$ by
\begin{equation}
\label{eq:iloczyn-rozlaczny-na-tensorach} (a_1\otimes \cdots \otimes
a_{\iloscklas})\bullet (b_1\otimes \cdots \otimes b_{\iloscklas})=
(a_1\bullet b_1)\otimes \cdots \otimes (a_{\iloscklas}\bullet
b_{\iloscklas}).
\end{equation}

Thus the map $\phi$ described in the above Lemma defines a
homomorphism of algebras $\phi:\A^{\otimes \iloscklas}\rightarrow
\C(G^q\rtimes S_q)$, where as the product in $\A^{\otimes
\iloscklas}$ we take the natural product.

We consider the map $\tilde{\phi}_{\zeta_i}:\A\rightarrow
\A^{\otimes \iloscklas}$ given by
$$\tilde{\phi}_{\zeta_i}(x)=\underbrace{1\otimes \cdots \otimes
1}_{i-1 \text{ times}}\otimes x \otimes \underbrace{1\otimes \cdots
\otimes 1}_{\iloscklas-i \text{ times}};$$ in this way
$\phi_{\zeta_i}=\phi\circ \tilde{\phi}_{\zeta_i}$. From the
following on we will not use any explicit numbering of the
representations in $\hat{G}$.

\subsection{Elements of the group algebra as random variables}
Let us fix some finite-dimensional representation $\rho_q$ of the
wreath product $G\wr S_q$. We can treat any commuting family of
elements of the group algebra $\C(G\wr S_q)$ as a family of random
variables equipped with the mean value given by the normalized
character:
\begin{equation}
\label{eq:wartoscoczekiwana} \E X:= \chi_{\rho_q}(X)= \tr \rho_q(X).
\end{equation}
It should be stressed that in the general case we treat elements of
$\C(G\wr S_q)$ as random variables only on a purely formal level; in
particular we do not treat them as functions on some Kolmogorov
probability space.

Lemma \ref{lem:glowny-homomorfizm} allows us to consider a
representation $\rho_q \circ \phi$ of $\A^{\otimes \iloscklas}$; for
simplicity we will denote this representation by the same symbol
$\rho_q$. In this way as random variables we may take the elements
of the algebra $\A^{\otimes \iloscklas}$. This algebra can be
equipped either with the natural product
\eqref{eq:iloczyn-naturalny-na-tensorach} and we denote the
resulting cumulants (called natural cumulants) by $k(X_1,\dots,X_n)$
or with the disjoint product $\bullet$
\eqref{eq:iloczyn-rozlaczny-na-tensorach} and we denote the
resulting cumulants (called disjoint cumulants) by
$k^{\bullet}(X_1,\dots,X_n)$.

\subsection{Canonical probability measure on Young diagrams associated to a representation}
\label{sub:Canonic-probabi-measure-Young-diagram-associa-represe}
There is a special case when it is possible to give a truly
probabilistic interpretation to the mean value
\eqref{eq:wartoscoczekiwana}: it is when for the family of random
variables we take the center of $\C(G\wr S_q)$ (respectively, the
algebra  $\A^{\otimes \iloscklas}$ equipped with the natural
product). The center of $\C(G\wr S_q)$ is isomorphic (via Fourier
transform) to the algebra of functions on irreducible
representations of $G\wr S_q$ and the expected value
\eqref{eq:wartoscoczekiwana} corresponds under this isomorphism to
the probability measure on irreducible representations (or,
equivalently, on functions $\Lambda:\hat{G}\rightarrow\mathbb{Y}$ as
described in Proposition \ref{prop:reprezentacje}) such that the
probability of a given irreducible representation is proportional to
the total dimension of all irreducible components of this type in
$\rho_q$.

\subsection{Generalized Young diagrams} \index{generalized Young
diagram} \index{Young diagram!generalized}
%

\begin{figure}[tb]
\includegraphics{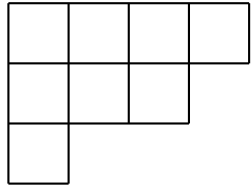}
\caption[Example of a Young diagram]{Young diagram associated with a
partition $8=4+3+1$.} \label{fig:young1}
\end{figure}

\begin{figure}[tb]
\includegraphics{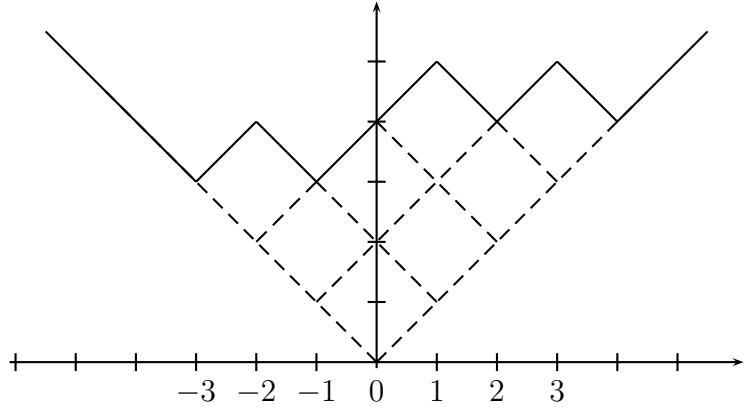}
\caption[Example of a generalized Young diagram]{Generalized Young
diagram associated with a partition $8=4+3+1$.} \label{fig:young2}
\end{figure}

Let $\lambda$ be a Young diagram. We assign to it a piecewise affine
function $\omega^\lambda:\R\rightarrow\R$ with slopes $\pm 1$, such
that $\omega^\lambda(x)=|x|$ for large $|x|$ as it can be seen on
the example from Figure \ref{fig:young2}. By comparing Figure
\ref{fig:young1} and Figure \ref{fig:young2} one can easily see that
the graph of $\omega^\lambda$ can be obtained from the graphical
representation of the Young diagram by an appropriate mirror image,
rotation and scaling by the factor $\sqrt{2}$. We call
$\omega^\lambda$ the generalized Young diagram associated with the
Young diagram $\lambda$ \cite{Kerov1993transition,
Kerov1998interlacing,Kerov1999differential}.


The class of generalized Young diagrams consists of all functions
$\omega:\R\rightarrow\R$ which are Lipschitz with constant $1$ and
such that $\omega(x)=|x|$ for large $|x|$ and of course not every
generalized Young diagram can be obtained by the above construction
from some Young diagram $\lambda$.

The setup of generalized Young diagrams is very useful in the study
of the asymptotic properties since it allows us to define easily
various notions of convergence of the Young diagram shapes.

\subsection{Functionals of the shape of Young diagrams}
The main result of this article is that the fluctuations of the
shape of some random Young diagrams converge (after some rescaling)
to a Gaussian distribution. Since the space of (generalized) Young
diagrams is infinite-dimensional therefore we need to be very
cautious when dealing with such statements. In fact, we will
consider a family of functionals on Young diagrams and we show that
the joint distribution of each finite set of these functionals
converges to the Gaussian distribution.

The functionals mentioned above are given as follows: for a Young
diagram $\lambda$ and the corresponding generalized Young diagram
$\omega$ we denote $\sigma(x)=\frac{\omega(x)-|x|}{2}$
\cite{Biane1998,IvanovOlshanski2002} and consider the family of maps
\begin{equation} \label{eq:p-tylda}
\tilde{p}_n(\lambda)= \int_{\R} x^n \sigma''(x) dx.
\end{equation}
Since $\sigma''$ makes sense as a distribution and $\sigma$ is
compactly supported hence the collection
$\big(\tilde{p}_n(\lambda)\big)_n$ determines the Young diagram
$\lambda$ uniquely.



\subsection{Transition measure of a Young diagram}
\index{transition measure of a Young diagram}
\label{subsec:transitionanalytic} To any generalized Young diagram
$\omega$ we can assign the unique probability measure $\mu^{\omega}$
on $\R$, called transition measure of $\omega$, the Cauchy transform
of which
\begin{equation}
\label{eq:Cauchy} G_{\mu^{\omega}}(z)= \int_{\R} \frac{1}{z-x} d
\mu^{\omega}(x)
\end{equation}
is given by
\begin{equation} \label{eq:definicja3} \log G_{\mu^{\omega}}(z)= -\frac{1}{2}
\int_{\R} \log(z-x)  \omega''(x) dx = -\frac{1}{2} \int_{\R}
\frac{1}{z-x} \omega'(x) dx
\end{equation}
%
%
for every $z\notin\R$. For a Young diagram $\lambda$ we will write
$\mu^{\lambda}$ as a short hand of $\mu^{\omega^{\lambda}}$. This
definition may look artificial but it turns out
\cite{Kerov1993transition,OkounkovVershik1996,Biane1998,
Okounkov2000randompermutations} that it is equivalent to natural
representation-theoretic definitions which arise by studying the
representation $\rho_q$ together with the inclusion $S_q\subset
S_{q+1}$.

For $p>0$ and a Young diagram $\lambda$ we consider the rescaled
(generalized) Young diagram $\omega^{p \lambda}$ given by $\omega^{p
\lambda}:x \mapsto p \omega^{\lambda}\big( \frac{x}{p} \big)$.
Informally speaking, the symbol $p \lambda$ corresponds to the shape
of the Young diagram $\lambda$ geometrically scaled by factor $p$
(in particular, if $\lambda$ has $q$ boxes then $p \lambda$ has $p^2
q$ boxes). It is easy to see that \eqref{eq:definicja3} implies that
the corresponding transition measure $\mu^{p \lambda}$ is a dilation
of $\mu^{\lambda}$:
\begin{equation}
\label{eq:skalowanie} \mu^{p \lambda}= D_p \mu^{\lambda}.
\end{equation}
This nice behavior of the transition measure with respect to
rescaling of Young diagrams makes it a perfect tool for the study of
the asymptotics of symmetric groups $S_q$ as $q\to\infty$.

\subsection{Free cumulants of the transition measure}
Cauchy transform of a compactly supported probability measure is
given at the neighborhood of infinity by a power series
$$ G_{\mu}(z)=\frac{1}{z}+\sum_{n\geq 2} M_n z^{-n-1},$$
where $M_n=\int_R x^n d\mu$ are the moments of the measure $\mu$. It
follows that on some neighborhood of infinity $G_{\mu}$ has a right
inverse $K_{\mu}$ with respect to the composition of power series
given by
$$ K_{\mu}(z)=\frac{1}{z}+\sum_{n\geq 1} R_n z^{n-1}$$
convergent on some neighborhood of $0$. The coefficients
$R_i=R_i(\mu)$ are called free cumulants of measure $\mu$. Free
cumulants appeared implicitly in Voiculescu's $R$--transform
\cite{VoiculescuAddition} and their combinatorial meaning was given
by Speicher \cite{Speicher1997}.

Free cumulants are homogenous in the sense that if $X$ is a random
variable and $c$ is some number then
$$ R_i(cX)=c^i R_i(X) $$
and for this reason they are very useful in the study of asymptotic
questions.

Each free cumulant $R_n$ is a polynomial in the moments
$M_1,M_2,\dots,M_n$ of the measure and each moment $M_n$ can be
expressed as a polynomial in the free cumulants $R_1,\dots,R_n$; in
other words the sequence of moments $M_1,M_2,\dots$ and the sequence
of free cumulants $R_1,R_2,\dots$ contain the same information about
the probability measure. The functionals of Young diagrams
considered in \eqref{eq:p-tylda} have a nice geometric
interpretation but they are not very convenient in actual
calculations. For this reason we will prefer to describe the shape
of a Young diagram by considering a family of functionals
\begin{equation}
\label{eq:wolnekumulanty} \lambda\mapsto R_n(\mu^{\lambda})
\end{equation}
given by the free cumulants of the transition measure. Equation
\eqref{eq:definicja3} shows that functionals $\tilde{p}_k$ from the
family \eqref{eq:p-tylda} can be expressed as polynomials in the
functionals from the family \eqref{eq:wolnekumulanty} and vice
versa.

Please note that the first two cumulants of a transition measure do
not carry any interesting information since
$$ R_1(\mu^{\lambda})=M_1(\mu^{\lambda})=0,$$
$$ R_2(\mu^{\lambda})=M_2(\mu^{\lambda})=q,$$
where $q$ denotes the number of the boxes of the Young diagram
$\lambda$.

Above we treated the free cumulant $R_i$ as a function on Young
diagrams, but it also can be viewed (via Fourier transform) as a
central element in $\C(S_q)$.

\section{The main result: representations with approximate factorization of characters}
\label{sec:main-result:-represe-with-approxi-factori-charact}

\subsection{The main theorems}
For a permutation $\sigma\in G^q\rtimes S_q$ we identify the coset
$\sigma G^q$ as an element of $(G^q\rtimes S_q)/G^q= S_q$. For a
permutation $\pi\in S_q$ we denote by $|\pi|$ the minimal number of
factors needed to write $\pi$ as a product of transpositions. For $l\leq q$ 
we treat $G^l \rtimes S_l$ as a subgroup of $G^q \rtimes S_q$.

The following theorem is the main result of this
article.
\begin{theoremanddefinition} \label{theo:main} For each
$q\geq 1$ let $\rho_q$ be a representation of $G\wr S_q$. We say
that the sequence $(\rho_q)$
has the character factorization property if it fulfills one (hence
all) of the following equivalent conditions:
\begin{itemize}

\item for any permutations $\sigma_1,\dots,\sigma_n\in G^l \rtimes
S_l$ with disjoint supports and such that
${\sigma_1}G^l,\dots,{\sigma_n}G^l\in S_l$ are cycles
\begin{equation}
\label{eq:faktoryzacjaHyper1} k(\sigma_1,\dots,\sigma_n)\ q^{
\frac{|{\sigma_1}G^l|+\cdots+|{\sigma_n}G^l|+2(n-1)}{2} }= O(1);
\end{equation}

\item for any integers $l_1,\dots,l_n\geq 1$ and any irreducible
representations $\zeta_1,\dots,\zeta_n\in\hat{G}$
\begin{equation}
\label{eq:faktoryzacjaHyper2}
k^{\bullet}\big(\tilde{\phi}_{\zeta_{1}}(\Sigma_{l_1}),\dots,\tilde{\phi}_{\zeta_{n}}(\Sigma_{l_n})\big)
\ q^{-\frac{l_1+\cdots+l_n-n+2}{2} } = O(1);
\end{equation}

\item for any integers $l_1,\dots,l_n\geq 1$ and any irreducible
representations $\zeta_1,\dots,\zeta_n\in\hat{G}$
\begin{equation}
\label{eq:faktoryzacjaHyper3}
k\big(\phi_{\zeta_1}(\Sigma_{l_1}),\dots,\phi_{\zeta_n}(\Sigma_{l_n})\big)
\ q^{-\frac{l_1+\cdots+l_n-n+2}{2} } = O(1);
\end{equation}

\item for any integers $l_1,\dots,l_n\geq 2$ and irreducible
representations $\zeta_1,\dots,\zeta_n\in\hat{G}$
\begin{equation}
\label{eq:faktoryzacjaHyper4} k\big(
\phi_{\zeta_1}(R_{l_1}),\dots,\phi_{\zeta_n} (R_{l_n})\big) \
q^{-\frac{l_1+\cdots+l_n-2(n-1)}{2}  }= O(1) .
\end{equation}
\end{itemize}
\end{theoremanddefinition}
Notice that in \eqref{eq:faktoryzacjaHyper2} we use disjoint cumulants which are
well defined in $\A^{\otimes \iloscklas}$ therefore the homomorphisms
$\tilde{\phi}_{\zeta_{j}}$ cannot be replaced by ${\phi}_{\zeta_{j}}$. 
On the other hand, in \eqref{eq:faktoryzacjaHyper3},
\eqref{eq:faktoryzacjaHyper4} we use the natural cumulants which are well
defined both in $\A^{\otimes \iloscklas}$ and in $\C(G^q \rtimes S_q)$ therefore
it does not matter if we use the homomorphism $\tilde{\phi}_{\zeta_{j}}$ or the
homomorphism ${\phi}_{\zeta_{j}}$.

The above theorem is a straightforward generalization of the analogous
result (Theorem and Definition 1 in
\cite{Sniady2005GaussuanFluctuationsofYoungdiagrams}) 
concerning representations of $\A$ to the case of
the representations of $\A^{\otimes \iloscklas}$ and for this reason
we skip its proof. Similarly, Theorem 3 in
\cite{Sniady2005GaussuanFluctuationsofYoungdiagrams} has the following analogue
for the wreath products:


\begin{theorem}
\label{theo:main2} Let $(\rho_q)$ has the character factorization
property. If the limit of one of the expressions
\eqref{eq:faktoryzacjaHyper1}--\eqref{eq:faktoryzacjaHyper4} exists
for $n\in\{1,2\}$ then the limits of all of the expressions
\eqref{eq:faktoryzacjaHyper1}--\eqref{eq:faktoryzacjaHyper4} exist
for $n\in\{1,2\}$.

These limits fulfill
\begin{multline}
\label{eq:granicamomentow} c_{\zeta,l+1}:= \lim_{q\to\infty}
\E\big(\phi_{\zeta }(\sigma,A_{\sigma})\big) q^{ \frac{l-1}{2} }=
\lim_{q\to\infty} \E\big(
\phi_{\zeta}(\Sigma_{l})\big) q^{-\frac{l+1}{2} } =\\
\lim_{q\to\infty} \E\big(\phi_{\zeta}(R_{l+1}) \big)
q^{-\frac{l+1}{2} },
\end{multline}
where $(\sigma,A_{\sigma})\in \partialS_q$ is a partial permutation
equal to a cycle of length $l$ and
\begin{multline}
\label{eq:granicamomentow2} \lim_{q\to\infty}
\Cov\big(\phi_{\zeta}(R_{l_1+1}),\phi_{\zeta}(R_{l_2+1}) \big)
q^{-\frac{l_1+l_2}{2}} = \\ \shoveright{\lim_{q\to\infty}\Cov\big(
\phi_{\zeta}(\Sigma_{l_1}),\phi_{\zeta}(\Sigma_{l_2})\big)
q^{-\frac{l_1+l_2}{2} }=}\\
\shoveleft{ \lim_{q\to\infty} \Cov^{\bullet} \big(
\tilde{\phi}_{\zeta}(\Sigma_{l_1}),\tilde{\phi}_{\zeta}(\Sigma_{l_2})\big)
q^{-
\frac{l_1+l_2}{2}  } }+ \\ \shoveright{ \sum_{r\geq 1}\sum_{\substack{a_1,\dots,a_r\geq 1\\
a_1+\dots+a_r=l_1}}  \sum_{\substack{b_1,\dots,b_r\geq 1\\
b_1+\dots+b_r=l_2}} \frac{l_1 l_2}{r} c_{\zeta,a_1+b_1} \cdots
c_{\zeta,a_r+b_r}=} \\
\shoveleft{ \lim_{q\to\infty}
\Cov\big(\phi_{\zeta}(\sigma_1,A_{\sigma_1}), \phi_{\zeta}
(\sigma_2,A_{\sigma_2}) \big)
 q^{\frac{l_1+l_2}{2} } -l_1 l_2 c_{\zeta,l_1+1} c_{\zeta,l_2+1}+} \\
\sum_{r\geq 1}\sum_{\substack{a_1,\dots,a_r\geq 1\\
a_1+\dots+a_r=l_1}}  \sum_{\substack{b_1,\dots,b_r\geq 1\\
b_1+\dots+b_r=l_2}} \frac{l_1 l_2}{r} c_{\zeta,a_1+b_1} \cdots
c_{\zeta,a_r+b_r},
\end{multline}
where $(\sigma_1,A_{\sigma_1}),(\sigma_2,A_{\sigma_2})\in
\partialS_q$ are disjoint cycles of length $l_1,l_2$, respectively,
and where the numbers $c_{\zeta,i}$ were defined in
\eqref{eq:granicamomentow}; also for any representations
$\zeta_1\neq \zeta_2\in \hat{G}$
\begin{multline}
\label{eq:granicamomentow2b} \lim_{q\to\infty} \Cov\big(
\phi_{\zeta_1}(R_{l_1+1}), \phi_{\zeta_2}(R_{l_2+1}) \big) q^{-\frac{l_1+l_2}{2}} =\\
\lim_{q\to\infty}\Cov \big( \phi_{\zeta_1}(\Sigma_{l_1}),
\phi_{\zeta_2}(\Sigma_{l_2}) \big)
q^{-\frac{l_1+l_2}{2} }=\\
\shoveleft{ \lim_{q\to\infty} \Cov^{\bullet}
\big( \tilde{\phi}_{\zeta_1}(\Sigma_{l_1}),\tilde{\phi}_{\zeta_2}(\Sigma_{l_2}) \big)
q^{- \frac{l_1+l_2}{2}  } }= \\
  \lim_{q\to\infty} \Cov\big(\phi_{\zeta_1}(\sigma_1,A_{\sigma_1}), \phi_{\zeta_2}
(\sigma_2,A_{\sigma_2}) \big) q^{\frac{l_1+l_2}{2} } -l_1 l_2
c_{\zeta_1,l_1+1} c_{\zeta_2,l_2+1}.
\end{multline}
\end{theorem}

By mimicking the line of the proof of Corollary 4 in
\cite{Sniady2005GaussuanFluctuationsofYoungdiagrams} we get the
following corollary.
\begin{corollary}
\label{coro:gaussian} Let $(\rho_q)$ be as in Theorem
\ref{theo:main2} and let $\Lambda:\hat{G}\rightarrow\mathbb{Y}$ be a
random function distributed according to the canonical probability
measure associated to $\rho_q$, as described in Section
\ref{sub:Canonic-probabi-measure-Young-diagram-associa-represe}.
\begin{enumerate}
\item
\label{item:gaba1}\emph{(Gaussian fluctuations of free cumulants)}
Then the joint distribution of the centered random variables
$$r_{\zeta,i}=q^{-\frac{i-1}{2}} (R_{\zeta,i} - \E R_{\zeta,i}) =
                q^{-\frac{i-1}{2}} \big(\phi_{\zeta}(R_{i}) - \E
\phi_{\zeta}(R_{i})\big)$$
converges to a Gaussian distribution in the weak topology of
probability measures, where $R_{\zeta,i}$ denotes the $i$-th free
cumulant of the transition measure $\mu^{\Lambda(\zeta)}$.

\item \label{item:gaba2} \emph{(Gaussian fluctuations of characters)}
Let $\sigma_i\in G\wr S_q$ be a permutation such that $\sigma_i
G^q\in S_q$ is a cycle. Then the joint distribution of the centered
random variables
\begin{equation}
\label{eq:gausowskiecharaktery} q^{ \frac{|\sigma_i G^q|+1}{2} }
\big(\chi_{\Lambda}(\sigma_i)- \E \chi_{\Lambda}(\sigma_i)\big)
\end{equation}
converges to a Gaussian distribution in the weak topology of
probability measures, where $\chi_{\Lambda} (x)=\tr \rho_{\Lambda}(x)$ denotes
the normalized character associated to the irreducible representation
$\rho_{\Lambda}$.
\item \label{item:gaba3}
\emph{(Gaussian fluctuations of the shape of the Young diagrams)}
Then the joint distribution of the centered random variables
$$\big( q^{-\frac{i-2}{2}} (\tilde{p}_{\zeta,i} - \E \tilde{p}_{\zeta,i})
\big)_{\zeta\in\hat{G}, i\geq 2} $$
converges to a Gaussian distribution in the weak topology of
probability measures, where
$\tilde{p}_{\zeta,i}=\tilde{p}_i\big(\Lambda(\zeta)\big)$ is the
functional of the shape of the Young diagram defined in
\eqref{eq:p-tylda}.
\end{enumerate}
\end{corollary}
%
%
%

\subsection{Examples}
All examples presented in this section not only have the character
factorization property but additionally are as in Theorem
\ref{theo:main2}, i.e.\ the limits of the means and of the
covariances exist.

\begin{example}
\label{example:1} Let $\rho:G\rightarrow \End(V)$ be a fixed
representation of $G$ and let $\rho_q=(V^{\otimes
q})\uparrow_{G^q}^{G^q\rtimes S_q}$ be the induced representation of
$G^q\rtimes S_q$. In particular, if $\rho$ is the left-regular
representation then $\rho_q$ is equal to the left regular
representation of the wreath product $G\wr S_q$.

It is easy to check that that for any permutations
$\sigma_1,\dots,\sigma_n\in G\wr S_q$ with disjoint supports
$$k(\sigma_1,\dots,\sigma_n) q^{
\frac{|\sigma_1 G|+\cdots+|\sigma_n G|+2(n-1)}{2} } = 0 \qquad
\text{if } n\neq 1 \text{ or } {\sigma_1}G\neq e. $$ It follows from
condition \eqref{eq:faktoryzacjaHyper1} that the sequence $(\rho_q)$
has the character factorization property and that the mean and the
covariance of the free cumulants are given by
\begin{equation} \label{eq:kerov1} \lim_{q\to\infty}
\E\big(\phi_{\zeta}(R_{l+1}) \big)
q^{-\frac{l+1}{2} } = \begin{cases} c_{\zeta} & \text{if } l=1 , \\
0 & \text{if } l\geq 2, \end{cases} \end{equation} where
$$c_{\zeta}= \frac{\text{(dimension of the representation $\zeta$)(multiplicity of $\zeta$ in $W$)} }%
{\text{(dimension of $W$)}} $$ and
\begin{multline} \label{eq:kerov2}
\lim_{q\to\infty} \Cov\big(\phi_{\zeta}(R_{l_1+1}),
\phi_{\zeta}(R_{\zeta,l_2+1}
) \big) q^{-\frac{l_1+l_2}{2}} = \\ \begin{cases} l c_{\zeta}^l & \text{if } l_1=l_2=l\geq 2, \\
c_{\zeta} ( 1-c_{\zeta}) & \text{if } l_1=l_2=1,
\\ 0 & \text{if } l_1\neq l_2,
\end{cases}
\end{multline}
and for $\zeta_1\neq \zeta_2$
\begin{multline}
\lim_{q\to\infty} \Cov\big(
\phi_{\zeta_1}(R_{l_1+1}),\phi_{\zeta_2}(R_{l_2+1}) \big)
q^{-\frac{l_1+l_2}{2}} = \\ \begin{cases}  - c_{\zeta_1}
c_{\zeta_2} & \text{if } l_1=l_2=1,  \\
0 & \text{otherwise.}
\end{cases}
\end{multline}
\end{example}

%

\begin{example}[Irreducible representations]
\label{example:3} Let $c>0$ be a constant and let $(\Lambda_q)$ be a
sequence of functions, $\Lambda_q:\hat{G}\rightarrow\mathbb{Y}$ as
in Proposition \ref{prop:reprezentacje}. We assume that
$\sum_{\zeta\in \hat{G}} |\Lambda_q(\zeta)|=q$ and that each diagram
$\Lambda_q(\zeta)$ has at most $c \sqrt{q}$ rows and columns.
Suppose that for each $\zeta\in\hat{G}$ the shapes of rescaled Young
diagrams $q^{-\frac{1}{2} } \Lambda_q(\zeta)$ converge to some
limit. The convergence of the shapes of Young diagrams implies
convergence of the free cumulants and it follows that the sequence
$(\rho_{\Lambda_q})$ of the corresponding irreducible
representations has the characters factorization property.

In this example the cumulants \eqref{eq:faktoryzacjaHyper3} and
\eqref{eq:faktoryzacjaHyper4} vanish for $n\geq 2$ since the Young
diagrams are non-random and the corresponding limits for $n=1$ are
determined by the limit shapes of the Young diagrams.
\end{example}

The above two examples are the building blocks from which one can
construct some more complex representations with the help of the
operations on representations presented below.

\begin{theorem}[Restriction of representations]
Suppose that the sequence of representations $(\rho_q)$ has the
character factorization property. Let a sequence of integers $(r_q)$
be given, such that $r_q\geq q$ and the limit ${p}=\lim_{q\to\infty}
\frac{q}{r_q} $ exists.

Let $\rho'_q$ denote the restriction of the representation
$\rho_{r_q}$ to the subgroup $G\wr S_q\subseteq G\wr S_{r_q}$. Then
the sequence $(\rho'_q)$ has the factorization property of
characters. The fluctuations of the free cumulants are determined by
\begin{equation}
c_{\zeta,l+1}':=\lim_{q\to\infty} \E(R'_{\zeta,l+1})
q^{-\frac{l+1}{2} }= {p}^{\frac{l-1}{2}} \lim_{q\to\infty}
\E(R_{\zeta,l+1}) q^{-\frac{l+1}{2} }= {p}^{\frac{l-1}{2}}
c_{\zeta,l+1}, \label{equ:restrictionA}
\end{equation}
\begin{multline} \lim_{q\to\infty} \Cov(R'_{\zeta,l_1+1},R'_{\zeta,l_2+1})
q^{-\frac{l_1+l_2}{2}}= \\ {p}^{\frac{l_1+l_2}{2}} \Bigg[
\lim_{q\to\infty} \Cov(R_{\zeta,l_1+1},R_{\zeta,l_2+1})
q^{-\frac{l_1+l_2}{2}} -
 l_1 l_2 c_{\zeta,l_1+1} c_{\zeta,l_2+1} \left( {p}^{-1}- 1\right) + \\
\sum_{r\geq 1}\sum_{\substack{a_1,\dots,a_r\geq 1\\
a_1+\dots+a_r=l_1}}  \sum_{\substack{b_1,\dots,b_r\geq 1\\
b_1+\dots+b_r=l_2}} \frac{l_1 l_2}{r} c_{\zeta,a_1+b_1} \cdots
c_{\zeta,a_r+b_r} \left(  {p}^{-r}-1 \right) \Bigg]
\label{mul:restrictionB}
\end{multline}
for all $\zeta\in\hat{G}$ and
\begin{multline} \lim_{q\to\infty} \Cov(R'_{\zeta_1,l_1+1},R'_{\zeta_2,l_2+1})
q^{-\frac{l_1+l_2}{2}}= \\ {p}^{\frac{l_1+l_2}{2}} \Bigg[
\lim_{q\to\infty} \Cov(R_{\zeta_1,l_1+1},R_{\zeta_2,l_2+1})
q^{-\frac{l_1+l_2}{2}} -
 l_1 l_2 c_{\zeta_1,l_1+1} c_{\zeta_2,l_2+1} \left( {p}^{-1}-
 1\right) \Bigg]
\label{mul:restrictionB2}
\end{multline}
for all $\zeta_1\neq \zeta_2$ and for all $l,l_1,l_2\geq 1$, where
the quantities $R'_{\zeta,i},c'_{\zeta,i}$ concern the
representations $(\rho'_q)$ while $R_{\zeta,i},c_{\zeta,i}$ concern
the representations $(\rho_q)$.

In particular, for ${p}=0$ we recover the fluctuations of the
induced representations from Example \ref{example:1}.
\end{theorem}
The proof is a straightforward application of Theorem and Definition
\ref{theo:main}; for details we refer to
Theorem 8 in \cite{Sniady2005GaussuanFluctuationsofYoungdiagrams}.

\begin{theorem}[Outer product of representations]
\label{theo:outerproduct} Suppose that for\/ $i\in\{1,2\}$ the
sequence of representations $(\rho^{(i)}_q)$ has the character
factorization property. Let sequences of positive integers
$r^{(i)}_q$ be given, such that $r^{(1)}_q+r^{(2)}_q=q$ and the
limits ${p}^{(i)}:=\lim_{q\to\infty} \frac{r^{(i)}}{q} $ exist.

Let $\rho'_q=\rho^{(1)}_{r^{(1)}_q} \circ \rho^{(2)}_{r^{(2)}_q}$
denote the outer product of representations. Then the sequence
$(\rho'_q)$ has the factorization property of characters with
\begin{equation}
\label{eq:outerA} c_{\zeta,l+1}':=\lim_{q\to\infty}
\E(R'_{\zeta,l+1}) q^{-\frac{l+1}{2} }=
\left({p}^{(1)}\right)^{\frac{l+1}{2}} c^{(1)}_{\zeta,l+1} +
\left({p}^{(2)}\right)^{\frac{l+1}{2}} c^{(2)}_{\zeta,l+1},
\end{equation}
and with an explicit (but involved) covariance of free cumulants.
The appropriate disjoint covariance is given by
\begin{multline*} \lim_{q\to\infty}
\Cov_{\rho'_q}^{\bullet}\big(\phi_{\zeta_1}(\Sigma_{l_1}),\phi_{\zeta_2}(\Sigma_{l_2})\big)
q^{-\frac{l_1+l_2}{2}}= \\( {p}^{(1)} )^{\frac{l_1+l_2}{2}}
\lim_{q\to\infty} \Cov_{\rho^{(1)}_{q}}^{\bullet}\big(
\phi_{\zeta_1}(\Sigma_{l_1}),\phi_{\zeta_2}(\Sigma_{l_2}) \big)
q^{-\frac{l_1+l_2}{2}}+ \\
({p}^{(2)} )^{\frac{l_1+l_2}{2}}  \lim_{q\to\infty}
\Cov_{\rho^{(2)}_{q}}^{\bullet}\big(\phi_{\zeta_1}(\Sigma_{l_1}),\phi_{\zeta_2}(\Sigma_{l_2})\big)
q^{-\frac{l_1+l_2}{2}}
\end{multline*}
for all $\zeta_1,\zeta_2\in\hat{G}$.
\end{theorem}
The proof of this result also follows closely the proof of the
analogous result for symmetric groups, namely Theorem 10 in
\cite{Sniady2005GaussuanFluctuationsofYoungdiagrams}.

\begin{theorem}[Induction of representations]
Suppose that the sequence of representations $(\rho_q)$ has
character factorization property. Let a sequence of integers $r_q$
be given, such that $r_q\leq q$ and the limit ${p}=\lim_{q\to\infty}
\frac{r_q}{q} $ exists.

Let $\rho'_q=\rho_{r_q}\big\uparrow_{G\wr S_{r_q}}^{G\wr S_q}$
denote the induced representation. Then the sequence $(\rho'_q)$ has
the characters factorization property  with
$$ c'_{\zeta,l+1}=\begin{cases} {p}^{\frac{l+1}{2}} c_{\zeta,l+1} & \text{for } l\geq 2, \\
p c_{\zeta,2} + (1-p) \frac{\text{(dimension of $\zeta$)}^2} {|G|} &
\text{for } l=1,
\end{cases}
$$ and with an explicit (but involved) covariance of free cumulants.
\end{theorem}
\begin{proof}
It is enough to adapt the proof of Theorem \ref{theo:outerproduct}.
\end{proof}

\begin{theorem}[Tensor product of representations]
Suppose that for $i\in\{1,2\}$ the sequence of representations
$(\rho^{(i)}_q)$ has character factorization property. Then the
tensor product $\rho'_q=\rho^{(1)}_q \otimes \rho^{(2)}_q$ has the
property of factorization of characters. Furthermore, the limit
distribution and the fluctuations are the same as for the tensor
representations \eqref{eq:kerov1} and \eqref{eq:kerov2} with
\begin{multline*}
c_{\zeta}=\lim_{q\to\infty} \frac{\left(\text{multiplicity of $\zeta$ in
$(\rho^{(1)}\otimes \rho^{(2)})\downarrow^{G\wr S_q}_G$}\right)
\left(\text{dimension of
$\zeta$}\right)}{\left(\text{dimension of $\rho^{(1)}\otimes
\rho^{(2)}$}\right)} .
\end{multline*}
\end{theorem}
Proof again follows the proof of the analogous result concerning symmetric
groups, namely Theorem 12 from
\cite{Sniady2005GaussuanFluctuationsofYoungdiagrams}.

\section{Acknowledgments}

I thank Akihito Hora for many discussions.

Research supported by the MNiSW research grant P03A 013 30, by the EU Research
Training Network \emph{QP-Applications}, contract HPRN-CT-2002-00279 and  by the
EC Marie Curie Host Fellowship for the Transfer of Knowledge \emph{Harmonic
Analysis, Nonlinear Analysis and Probability}, contract MTKD-CT-2004-013389.

\bibliographystyle{alpha}
\bibliography{biblio}

\end{document}